\documentclass[10pt,twocolumn]{article}

\usepackage{amsmath}
\usepackage{amssymb}
\usepackage{color}
\usepackage{graphicx}
\usepackage[permil]{overpic}
\usepackage[a4paper, total={7.2in, 10.2in}]{geometry}

\newcommand{\comenta}[1]{}

\newcommand{\az}{\textcolor{black}}

\newcommand{\daniel}{\textcolor{black}}

\title{\LARGE \bf
Distribution of a Markov chain in reverse-time with cluster observations in the extremes of a finite time window}
\author{Daniel A. Gutierrez-Pachas, Eduardo F. Costa and Alessandro N. Vargas % <-this % stops a space
  \thanks{Daniel A. Gutierrez-Pachas is with the Department of Computer Science. Universidad Católica San Pablo, Arequipa, Peru
        {\tt\small daniel.gutierrez.inv@ucsp.edu.pe}.
Eduardo F. Costa is with the Instituto de Ciências Matemáticas e de Computação, Universidade de São Paulo, São Carlos, Brazil
         {\tt\small efcosta@icmc.usp.br}.
 Alessandro N. Vargas is with the Departamento de Engenharia Elétrica, Universidade Tecnológica Federal do Paraná, Brazil
         {\tt\small avargas@utfpr.edu.br}. This work was supported by
CNPq Grant 310877/2017-2, CNPq-Universal 421486/2016-3,
FAPESP 2017/20934-9 and FAPESP-CEPID 2013/07375-0.
}
 }

\date{}

\begin{document}

\maketitle
\thispagestyle{empty}
\pagestyle{empty}

\begin{abstract}
In this brief note, we find formulas for the distribution and the transition probability matrices of a stochastic process described as a time-reversion in a finite time window of a Markov chain, with cluster observation of the Markov state in the extremes of that window.  \end{abstract}

%\section{Basic notation and information structure}\label{sec-nota}
\section*{Problem formulation and solution}
Consider a Markov chain $\{\eta(k),\,\,k=0,1,\ldots \}$ 
%defined in a fixed  probability space $(\Omega,\mathcal{F},\text{Pr})$,  %\cite{Book-measure-theory} 
with distribution and transition probabilities given respectively by 
\begin{equation}\nonumber
\begin{aligned}
s_{ij}(k)&=\text{Pr}(\eta(k+1)=j|\eta(k)=i), 
\\ 
\upsilon_i(k)&=\text{Pr}(\eta(k)=i),\;\;i,j\in\mathbb{N}, k\geq 0,
\end{aligned}    
\end{equation}
where $\mathbb{N}=\{1,\ldots,N\}$ is the state space. 
For Markov chains and basic probability rules employed in this note, the reader is referred to 
\cite{Cinlar75,Papoulis}.
Consider the process $\{\theta(k), k=0,1,\ldots,\ell\}$, defined as follows, based on the  Markov chain in a finite time window,
$$\theta(k)=\eta(\ell-k), \;\; \,k=0,1,\ldots,\ell.$$ 
Assume that a cluster observation is available 
in the form $\{\theta(0)\in \mathbb{C}_0\}$, $\mathbb{C}_0\subseteq \mathbb{N}$, as well as a cluster observation with anticipation 
$\{\theta(\ell)\in \mathbb{C}_\ell$\}, 
$\mathbb{C}_\ell\subseteq\mathbb{N}$.
$E_{\text{o}}\in\mathcal{F}$ stands for the associated $\sigma$-algebra. 

The adopted process $\theta$ and information structure may appear in real-world problems, including problems involving a first-in last-out queue. For example, suppose the Markov chain characterise the presence of certain irregularities in $\ell$ items arriving at a store, and that the first item in goes through an inspection of these irregularities, so that we know the state of this item (when it enters the queue); following the queue rule, and counting the items when leaving the queue, we have information ``with anticipation'' on the state of the $\ell$-th item. Note, via this example, that observation with anticipation does not necessarily require prediction of future events.

In this note we seek for formulas 
for the conditional transition probabilities and distributions 
defined in (1), based on the Markov chain parameters $\upsilon_i$ and \daniel{$s_{ij}(k)$} and the cluster observations. We shall need some additional notation.
For each $i,j\in\mathbb{N}$ and  $k=0,1,\ldots,\ell-1$, we define
\begin{equation}\label{eq-prob-theta}
\begin{aligned}
p_{ij}&=\text{Pr}(\theta(k)=j\,|\,\theta(k+1)=i\,,\,E_{\text{o}}),
\\
\pi_i(k)&=\text{Pr}(\theta(k)=i)|E_{\text{o}}), 
\;\;i,j\in\mathbb{N}, 0\leq k\leq \ell.
\end{aligned}
\end{equation}
We write $\mathbf{P}(k)$ to represent a matrix of dimension $N$ by $N$, whose components are $p_{ij}(k)$, therefore satisfying the Chapman-Kolmogorov equation 
\begin{equation}\label{eq-chapman}
\pi(k+1)=\mathbf{P}(k)\pi(k).
\end{equation}

%-----------------------------------------------------------
\noindent{W}e denote, for $i\in\mathbb{N}$, and $k=1,\ldots,\ell-1$, 
$${e}= \sum_{r\in\mathbb{C}_\ell}\sum_{j\in\mathbb{C}_{0}}[\mathbf{S}^{\ell}]_{rj}\upsilon_r(\az{0}),\quad {g}_i(k)=\sum_{r\in\mathbb{C}_{0}}[\mathbf{S}^{k+1}]_{ir}.$$
%\az{and $f_i(k)$ is the $i-th$ element of the vector $\upsilon(0)\mathbf{S}^{\ell-k-1}$.}

%-----------------------------------------------------------
%-----------------------------------------------------------
\noindent\textbf{Lemma 1.} 
\textit{The following is valid for $1\leq k\leq \ell-2$ and $i,j\in\mathbb{N}$; if \az{$i$ is such that $\upsilon_i(\ell-k-1)>0$}, then
\begin{equation}\label{eq-prob_lema3}
    p_{ij}(k)=\begin{cases}
    (g_i(k))^{-1}\sum_{r\in\mathbb{C}_0}[\mathbf{S}^{k}]_{jr}s_{ij}(k),&\, g_i(k)\neq 0,\\
    \quad\quad\quad \text{arbitrary},&\, {otherwise};
    \end{cases}
\end{equation}
\az{if  $\upsilon_i(\ell-k-1)=0$, then  $p_{ij}(k)$ is arbitrary.} 
\eqref{eq-prob_lema3} is also valid for $k=\ell-1$, $i\in\mathbb{C}_{\ell}$ and $j\in\mathbb{N}$, as well as for $k=0$, $i\in\mathbb{N}$ and $j\in\mathbb{C}_{0}$. %
The remaining cases regarding $p_{ij}(k)$ are: 
$p_{ij}(\ell-1)$ \az{is arbitrary} for $i\notin\mathbb{C}_{\ell}$ and $j\in\mathbb{N}$; 
\az{if $i\in\mathbb{C}_\ell$ and there is 
no Markov state $\eta(\ell)\in\mathbb{C}_{0}$ that can be reached from $i$ 
(in the sense that $\text{Pr}(\eta(\ell)\in\mathbb{C}_{0}\,|\,\eta(0)=i)=0$) then $p_{ij}(\ell-1)$ is  arbitrary;} 
$p_{ij}(0)$ is arbitrary for $i\in\mathbb{N}$ and $j\notin\mathbb{C}_{0}$. Regarding $\pi_{i}(\ell)$, for $i\in\mathbb{C}_{\ell}$, 
\begin{equation}\label{eq-pi-terminal}
    \pi_{i}(\ell)=\begin{cases}
    e^{-1}\displaystyle\sum_{j\in\mathbb{C}_{0}} [\mathbf{S}^{\ell}]_{ij}\upsilon_i(0),&\, e\neq 0,\\
    \quad\quad\quad \text{arbitrary},&\, {otherwise}.
    \end{cases}
\end{equation}
For $i\notin\mathbb{C}_{\ell}$, $\pi_{i}(\ell)=0$. 
Finally, $\pi(k)$, $1\leq k\leq \ell-1$, is given by \eqref{eq-chapman} and the formulas above.
}

\textbf{Proof:} for $i\in\mathbb{N}$ by definition we have 
\begin{equation}\nonumber
\begin{aligned}
 \pi_i(\ell)&=\text{Pr}(\theta(\ell)=i\,|\,\eta(0)\in\mathbb{C}_\ell,\,\eta(\ell)\in\mathbb{C}_0)\\
 &=\text{Pr}(\eta(0)=i\,|\,\eta(0)\in\mathbb{C}_\ell,\,\eta(\ell)\in\mathbb{C}_0).  
  \end{aligned}
\end{equation}
Of course,  
$\pi_i(\ell)=0$  whenever $i\notin\mathbb{C}_\ell$, as it is given that 
$\theta(\ell)=\eta(0)$ is in $\mathbb{C}_\ell$. If $i\in\mathbb{C}_\ell$, and assuming that 
$\text{Pr}\big(\eta(0)\in\mathbb{C}_\ell\,,\,\eta(\ell)\in\mathbb{C}_0\big)>0$, we may write 
\begin{equation*}
\begin{aligned}
\pi_i(\ell)&=\text{Pr}\big(\eta(0)=i\,|\,\eta(0)\in\mathbb{C}_\ell\,,\,\eta(\ell)\in\mathbb{C}_0\big)\\
    &=\frac{\text{Pr}\big(\eta(\ell)\in\mathbb{C}_0\,\,,\,\,\eta(0)=i\big)}{\text{Pr}\big(\eta(\ell)\in\mathbb{C}_0\,\,,\,\,\eta(0)\in\mathbb{C}_\ell\big)}\\
    &=\frac{\displaystyle\sum_{j\in\mathbb{C}_{0}}\text{Pr}\big(\eta(\ell)=j\,,\,\eta(0)=i\big)}{\displaystyle\sum_{r\in\mathbb{C}_{\ell}}\sum_{j\in\mathbb{C}_0}\text{Pr}\big(\eta(\ell)=j\,,\,\eta(0)=r\big)}.
\end{aligned}
\end{equation*}
Note that
\begin{equation}\nonumber
\begin{aligned}
\text{Pr}(\eta(\ell)=j,\,\eta(0)=r)&=\text{Pr}(\eta(\ell)=j|\eta(0)=r)\\
&\cdot\text{Pr}(\eta(0)=r)=[\mathbf{S}^{\ell}]_{rj}\upsilon_r(0), 
\end{aligned}
\end{equation}
and substituting in the above
$$\begin{aligned}
\pi_i(\ell)&=\displaystyle\sum_{j\in\mathbb{C}_{0}}[\mathbf{S}^{\ell}]_{i,j}\upsilon_i(0)\left(\displaystyle\sum_{r\in\mathbb{C}_\ell}\sum_{j\in\mathbb{C}_{0}}[\mathbf{S}^{\ell}]_{r,j}\upsilon_r(0)\right)^{-1}
\\& 
=e^{-1}\displaystyle\sum_{j\in\mathbb{C}_{0}} [\mathbf{S}^{\ell}]_{ij}\upsilon_i(0)
\end{aligned}$$
If $i\in\mathbb{C}_\ell$ and 
$\text{Pr}\big(\eta(0)\in\mathbb{C}_\ell\,,\,\eta(\ell)\in\mathbb{C}_0\big)=0.$ 
(in which case, the inverse in the above equation does 
not exist), then $\pi_i(\ell)$ 
%&=\text{Pr}\big(\eta(0)=i\,|\,\eta(0)\in\mathbb{C}_\ell\,,\,\eta(\ell)\in\mathbb{C}_0\big)$ 
is conditioned on an event of probability zero 
therefore $\pi_i(\ell)$ is arbitrary,
thus completing the demonstration of \eqref{eq-pi-terminal}.
%
%whenever $\displaystyle\sum_{r\in\mathbb{C}_\ell}\sum_{j\in\mathbb{C}_{0}}[\mathbf{S}^{\ell}]_{r,j}\upsilon_r(0)\neq 0$, otherwise $\pi_i(\ell)=0$. 
%
Now we turn our attention to $\mathbf{P}$.
\begin{equation}\nonumber
\begin{aligned}
p_{ij}(\ell-1)&=\text{Pr}(\theta(\ell-1)=j|\theta(\ell)=i,\eta(0)\in\mathbb{C}_\ell,\eta(\ell)\in\mathbb{C}_0)\\
&=\text{Pr}(\eta(1)=j|\eta(0)=i,\eta(0)\in\mathbb{C}_\ell,\eta(\ell)\in\mathbb{C}_0).
\end{aligned}
\end{equation}
It is clear that, if $i\notin\mathbb{C}_\ell$, then $p_{ij}(\ell-1)$ 
is conditioned  on an empty set, therefore an event of probability zero, making 
$p_{ij}(\ell-1)$ arbitrary, for any $j\in\mathbb{N}$. 
If $i\in\mathbb{C}_\ell$ and the event $\{\eta(\ell)\in\mathbb{C}_{0}\,|\,\eta(0)=i\}$ 
is not of probability zero, then using the total probability law we write
%
% If $\displaystyle\sum_{r\in\mathbb{C}_{\ell}}[\mathbf{S}^{\ell}]_{ir}\neq 0$, we obtain
\begin{equation}\nonumber
\begin{aligned}
p_{ij}(\ell-1)
&=\frac{\text{Pr}(\eta(\ell)\in\mathbb{C}_{0}\,|\,\eta(1)=j,\eta(0)=i)\cdot s_{ij}(k)} {\text{Pr}(\eta(\ell)\in\mathbb{C}_{0}\,|\,\eta(0)=i)}\\
&=\frac{\text{Pr}(\eta(\ell)\in\mathbb{C}_{0}\,|\,\eta(1)=j)\cdot  s_{ij}(k)}
 {\text{Pr}(\eta(\ell)\in\mathbb{C}_{0}\,|\,\eta(0)=i)}\\
&=\displaystyle\sum_{r\in\mathbb{C}_{\ell}}[\mathbf{S}^{\ell-1}]_{jr} s_{ij}(k) \left(\displaystyle\sum_{r\in\mathbb{C}_{\ell}}[\mathbf{S}^{\ell}]_{ir}\right)^{-1},  j\in\mathbb{N},
\end{aligned}
\end{equation}
where the second inequality is due to the Markov property. 
If $i\in\mathbb{C}_\ell$ and the event $\{\eta(\ell)\in\mathbb{C}_{0}\,|\,\eta(0)=i\}$ 
is an empty set, then $p_{ij}(\ell-1)$ is conditional on an event of zero probability, hence it is arbitrary. 
Regarding $p_{ij}(k)$ with $1\leq k\leq \ell-2$ and $i,j\in\mathbb{N}$, 
denoting 
$A = \{\eta(\ell-k-1)=i\}$ for a better visual diagramming of 
the next equation, 
assuming $\text{Pr}(\eta(\ell)\in\mathbb{C}_{0},A)>0$ and using the Markov property, we write 
\begin{equation}\nonumber
\begin{aligned}
& p_{ij}(k)=\text{Pr}(\eta(\ell-k)=j|A,\eta(0)\in\mathbb{C}_\ell,\eta(\ell)\in\mathbb{C}_0)\\
& \;\;=\text{Pr}(\eta(\ell-k)=j|A,\eta(\ell)\in\mathbb{C}_{0})\\
& \;\;=\frac{\text{Pr}(\eta(\ell-k)=j,\eta(\ell)\in\mathbb{C}_{0},A)}{\text{Pr}(\eta(\ell)\in\mathbb{C}_{0},A)}\\
& \;\;=\frac{\text{Pr}(\eta(\ell)\in\mathbb{C}_{0}|\eta(\ell-k)=j,A)\text{Pr}(\eta(\ell-k)=j|A)P(A)} {\text{Pr}(\eta(\ell)\in\mathbb{C}_{0}|A)P(A)}\\
& \;\;=\frac{\text{Pr}(\eta(\ell)\in\mathbb{C}_{0}|\eta(\ell-k)=j)\text{Pr}(\eta(\ell-k)=j|A)} {\text{Pr}(\eta(\ell)\in\mathbb{C}_{0}|A)},
\end{aligned}
\end{equation}
yielding
$$\begin{aligned}
p_{ij}(k)&=\sum_{r\in\mathbb{C}_{0}}[\mathbf{S}^{k}]_{jr} s_{ij}(k) \left(\displaystyle\sum_{r\in\mathbb{C}_{0}}[\mathbf{S}^{k+1}]_{ir}\right)^{-1}
\\&
=(g_i(k))^{-1}\sum_{r\in\mathbb{C}_0}[\mathbf{S}^{k}]_{jr}s_{ij}(k).
\end{aligned}
$$
\az{Note that the requirement $\text{Pr}(\eta(\ell)\in\mathbb{C}_{0},A)>0$ is equivalent to say: (i) considering the Markov chain $\{\eta_k, k\geq 0\}$,  
the set $\mathbb{C}_{0}$ is reachable from state $i$ in $k+1$ steps 
(in which case $g_i(k)=0$), and (ii) $\text{Pr}(\eta(\ell-k-1)=i)=\upsilon_i(\ell-k-1)>0$. If any of (i) or (ii) is false, or both are false, then then $p_{ij}(k)$  is conditional on an event of probability zero, 
hence it is \az{arbitrary}.}
It only remains to find the formula for $p_{ij}(0)$.
\begin{equation}\nonumber
    \begin{aligned}
p_{ij}(0)&=\text{Pr}(\theta(0)=j|\theta(1)=i,\eta(0)\in\mathbb{C}_\ell,\eta(\ell)\in\mathbb{C}_0)\\
&=\text{Pr}(\eta(\ell)=j|\eta(\ell-1)=i,\eta(0)\in\mathbb{C}_\ell,\eta(\ell)\in\mathbb{C}_0).
    \end{aligned}
\end{equation}
If  $j\notin\mathbb{C}_{0}$ 
then $p_{ij}(0)=0$ because it is given 
that $\eta(0)\in\mathbb{C}_{0}$; otherwise, we have
\begin{equation}\nonumber
\begin{aligned}
p_{ij}(0)&=\text{Pr}\big(\eta(\ell)=j\,|\,\eta(\ell-1)=i\,,\,\eta(\ell)\in\mathbb{C}_0\big)\\
\end{aligned}
\end{equation}
so that, when the conditional event is not of probability zero, 
\begin{equation}\nonumber
\begin{aligned}
%p_{ij}(0)&=\text{Pr}\big(\eta(\ell)=j\,|\,\eta(\ell-1)=i\,,\,\eta(\ell)\in\mathbb{C}_0\big)\\
p_{ij}(0)&=\frac{\text{Pr}\big(\eta(\ell)=j,\eta(\ell)\in\mathbb{C}_0\,|\,\eta(\ell-1)=i\big)}{\text{Pr}\big(\eta(\ell)\in\mathbb{C}_0\,|\,\eta(\ell-1)=i\big)}\\
&=\frac{\text{Pr}\big(\eta(\ell)=j\,|\,\eta(\ell-1)=i\big)}{\text{Pr}\big(\eta(\ell)\in\mathbb{C}_0\,|\,\eta(\ell-1)=i\big)}\\
&=s_{ij}(k) \left(\displaystyle\sum_{r\in\mathbb{C}_{0}}s_{ir}\right)^{-1},\forall\,i\in\mathbb{N}\,\,\text{and}\,\,j\in\mathbb{C}_{0}, 
\end{aligned}
\end{equation}
and, when the conditional event $\{\eta(\ell-1)=i,\,\eta(\ell)\in\mathbb{C}_0\}$ is of probability zero 
(in which case the inverse in the above equation does not exist), then $p_{ij}(0)$ is \az{arbitrary}.
  
Finally, $\pi(k)$, $1\leq k\leq \ell-1$, is given by \eqref{eq-chapman} and the formulas above.
\hfill Q.E.D.

\medskip
{\it Remark 1.} All arbitrary values in Lemma 1 can be set to zero. This is a suitable choice 
in some cases, as in \cite{Daniel-duality,PachasSubm}, where a Riccati-like equation is computed for every $i,k$ such that $\pi_i(k)>0$, so that, choosing $\pi_i(k)=0$ avoids unnecessary computations. 

\medskip
{\it Remark 2.} In view of \eqref{eq-prob_lema3},
the transition probabilities of the process $\theta$ depend on $k$ even if the Markov chain is time-homogeneous. 
For $\mathbf{P}$ to be irrespective of $k$, it would be necessary that the Markov chain is time-homogeneous \emph{and} there is no observation of $\theta(0)$,  that is, $\mathbb{C}_0=\mathbb{N}$; in this case, 
\eqref{eq-prob_lema3} reduces to $p_{ij}=s_{ij}$.

%\bibliographystyle{IEEEtran}
%\bibliography{mjls}

\begin{thebibliography}{1}
\providecommand{\url}[1]{#1}
\csname url@samestyle\endcsname
\providecommand{\newblock}{\relax}
\providecommand{\bibinfo}[2]{#2}
\providecommand{\BIBentrySTDinterwordspacing}{\spaceskip=0pt\relax}
\providecommand{\BIBentryALTinterwordstretchfactor}{4}
\providecommand{\BIBentryALTinterwordspacing}{\spaceskip=\fontdimen2\font plus
\BIBentryALTinterwordstretchfactor\fontdimen3\font minus
  \fontdimen4\font\relax}
\providecommand{\BIBforeignlanguage}[2]{{%
\expandafter\ifx\csname l@#1\endcsname\relax
\typeout{** WARNING: IEEEtran.bst: No hyphenation pattern has been}%
\typeout{** loaded for the language `#1'. Using the pattern for}%
\typeout{** the default language instead.}%
\else
\language=\csname l@#1\endcsname
\fi
#2}}
\providecommand{\BIBdecl}{\relax}
\BIBdecl

\bibitem{Cinlar75}
E.~Cinlar, \emph{Introduction to Stochastic Processes}.\hskip 1em plus 0.5em
  minus 0.4em\relax Prentice-Hall, 1975.

\bibitem{Papoulis}
A.~Papoulis, \emph{Probability, Random Variables, and Stochastic
  Processes}.\hskip 1em plus 0.5em minus 0.4em\relax New York: McGraw-Hill,
  1984.

\bibitem{Daniel-duality}
D.~A. Gutierrez-Pachas and E.~F. Costa, ``On the linear quadratic problem for
  systems with time-reversed {M}arkov jump parameters and the duality with
  filtering of {M}arkov jump linear systems,'' \emph{IEEE Transactions on
  Automatic Control}, vol.~63, no.~9, pp. 3040--3045, 2018.

\bibitem{PachasSubm}
D.~A. Gutierrez-Pachas, E.~F. Costa, and A.~N. Vargas, ``Linear quadratic
  control problem of systems with {M}arkov jumps in reverse time and
  observation with anticipation of the jumps,'' \emph{Submitted to the IEEE
  Control Systems Letters}.

\end{thebibliography}

\end{document}